\documentclass[12pt,a4paper]{article}
\usepackage{amssymb,latexsym,a4wide}
\def\exp{\mathop{\rm exp}\nolimits}

\def\trdeg{\mathop{\rm trdeg}\nolimits}

\def\ZZ{{\mathbb Z}}
\def\RR{{\mathbb R}}

\def\bull{\vrule height .9ex width .8ex depth -.1ex }

\newtheorem{formula}{}[section]
\newtheorem{proposition}[formula]{Proposition}
\newtheorem{definition}[formula]{Definition}
\newtheorem{corollary}[formula]{Corollary}
\newtheorem{remark}[formula]{Remark}
\newtheorem{lemma}[formula]{Lemma}
\newtheorem{theorem}[formula]{Theorem}

\def\thrm{\begin{theorem}}
\def\thrml#1{\begin{theorem}\label{#1}}
\def\ethrm{\end{theorem}}
\def\rmrk{\begin{remark}}
\def\rmrkl#1{\begin{remark}\label{#1}}
\def\ermrk{\end{remark}}
\def\dfntn{\begin{definition}}
\def\dfntnl#1{\begin{definition}\label{#1}}
\def\edfntn{\end{definition}}
\def\nmrt{\begin{enumerate}}
\def\enmrt{\end{enumerate}}

\def\qtn{\begin{equation}}
\def\qtnl#1{\begin{equation}\label{#1}}
\def\eqtn{\end{equation}}
\def\lmm{\begin{lemma}}
\def\lmml#1{\begin{lemma}\label{#1}}
\def\elmm{\end{lemma}}
\def\crllr{\begin{corollary}}
\def\crllrl#1{\begin{corollary}\label{#1}}
\def\ecrllr{\end{corollary}}
\begin{document}
\title{On a tropical dual Nullstellensatz}


\author{
Dima Grigoriev \\[-1pt]
\small CNRS, Math\'ematiques, Universit\'e de Lille \\[-3pt]
\small Villeneuve d'Ascq, 59655, France\\[-3pt]
{\tt \small Dmitry.Grigoryev@math.univ-lille1.fr}\\[-3pt]
\small http://en.wikipedia.org/wiki/Dima\_Grigoriev
}
\date{}
\maketitle

\begin{abstract}
Since a tropical Nullstellensatz fails even for tropical univariate
polynomials we study a conjecture on a tropical {\it dual}
Nullstellensatz for tropical polynomial systems in terms of
solvability of a tropical linear system with the Cayley matrix
associated to the tropical polynomial system. The conjecture on a
tropical effective dual Nullstellensatz is proved for tropical
univariate polynomials.
\end{abstract}


\section*{Introduction}

Let $T$ be a tropical semi-ring with operations $\oplus,\, \otimes$
(see e.~g. \cite{Butkovic}, \cite{Sturmfels}, \cite{Mikhalkin},
\cite{Theobald}). Typically $\oplus=\min,\, \otimes=+$. Examples of
$T$ are $\ZZ$ and $\ZZ_{\infty}=\ZZ\cup \{\infty\}$. A {\it tropical
monomial} has a form $Q=a\otimes X_1^{\otimes i_1}\otimes \cdots
\otimes X_n^{\otimes i_n}, \, a\in T$. The {\it tropical degree}
$\trdeg (Q):= i_1+\cdots+i_n$. From the point of view of the
classical algebra a tropical monomial is a linear function. A point
$x=(x_1,\dots,x_n)\in T^n$ (with some of $x_i \neq \infty$) is a
{\it tropical zero} of a {\it tropical polynomial} $f=\bigoplus _l
Q_l$ if the minimum $\min_l \{Q_l(x)\}$ is attained for at least two
different tropical monomials $Q_l$.

We study the issue of a tropical Nullstellensatz. Its direct
formulation fails even for tropical univariate polynomials: for
example, two tropical polynomials $X\oplus 0,\, X\oplus 1$ have no
common tropical zero, while the generated by them tropical ideal
does not contain $1$ or any other tropical monomial. That is why we
consider a {\it tropical "dual" Nullstellensatz}.

One can treat the (customary) Hilbert's Nullstellensatz as a
reduction of solvability of a polynomial system to solvability of a
suitable linear system. Namely, solvability of a polynomial system
is equivalent to that the Cayley matrix $C$ associated to the system
does not contain the vector $(1,0,\dots,0)$ in the linear hull of
its rows. In its turn it is equivalent to that the linear system
$C\cdot (a_0,a_1,\dots)=0$ has a solution with $a_0\neq 0$ (cf.
Section~\ref{section1}). The latter rephrasing of the
Nullstellensatz we call the "dual" Nullstellensatz. It holds also
for the {\it infinite} matrix $C$ (we call it the infinite "dual"
Nullstellensatz) unlike the customary Nullstellensatz, and it holds
for a finite submatrix of $C$ with the size depending on $n$ and on
the degrees of the polynomials in the system (we call it the {\it
effective} "dual" Nullstellensatz).

In Section~\ref{section2} we formulate the conjecture on a tropical
"dual" Nullstellensatz. In Section~\ref{section3} we give a
rephrasing of the conjecture in terms of the combinatorial convex
geometry. Finally, in Section~\ref{section4} we prove the tropical
effective "dual" Nullstellensatz for univariate polynomials ($n=1$).

Observe that the latter result in case of a system of two tropical
polynomials $f_1,\, f_2$ follows from the approach of \cite{Tabera}
which relies on the (classical) resultant of a pair of (classical)
polynomials, so it fails for overdetermined systems in the tropical
setting. We mention also that the problem of solvability of tropical
polynomial systems is $NP$-complete even for tropical quadratic
polynomials \cite{Theobald}.

Solvability of tropical linear systems belongs to the complexity
class $NP\cap co-NP$ \cite{Akian}, \cite{G}. In \cite{Akian},
\cite{G} two different algorithms for solving tropical linear
systems were designed with the similar complexity bounds polynomial
in $s.\, M$, where $s$ is the size of the tropical linear system
(so, of its matrix) and $M$ majorates the absolute values of the
finite (integer) coefficients of the system. We note that the
algorithm from \cite{G} possesses an extra feature that it has also
a complexity bound polynomial in $\exp(s),\, \log M$. The open
question is whether it runs in fact, within complexity polynomial in
$s, \, \log M$ (which would provide a polynomial complexity for the
problem of solvability of tropical linear systems)?

In addition, the algorithm from \cite{G} entails as a by-product the
equivalence of solvability of a tropical linear system with the
degeneration of its {\it tropical rank} and simultaneously with the
degeneration of its {\it Kapranov rank}. The latter for systems with
{\it finite} coefficients (say, from $\ZZ$) was shown in
\cite{Sturmfels}, also a part of this equivalence just for the
tropical rank follows from \cite{Rowen}.

Besides, we mention that in \cite{Izhakian} the tropical (customary)
Nullstellensatz was established for an introduced there a "ghost"
tropical semi-ring. In \cite{Shustin} the radical of a tropical
ideal was explicitly described.

\section{"Dual" Nullstellensatz}\label{section1}

Let $F_1,\dots,F_s\in K[X_1,\dots,X_n]$ be polynomials over an
algebraically closed field $K$. Denote by $C:=C(F_1,\dots,F_s)$ the
(infinite size) Cayley matrix over $K$ consisting of the
coefficients of $F_1,\dots,F_s$. The columns of $C$ correspond to
all the monomials $X^I:=X_1^{i_1}\cdots X_n^{i_n},\,
I=(i_1,\dots,i_n)$, and the rows of $C$ correspond to all the
polynomials of the form $X^I\cdot F_j, 1\le j\le s$. Let the first
column of $C$ correspond to the monomial $X^0=1$. For an integer $N$
denote by $C_N$ the (finite size) submatrix of $C$ formed by the
rows $X^I\cdot F_j,\, 1\le j\le s$ with the degrees $\deg
X^I=i_1+\cdots+i_n\le N$ and the corresponding columns which contain
a non-zero entry in at least one of these rows.

Nullstellensatz states that a polynomial system
\begin{eqnarray}\label{1}
F_1=\cdots=F_s=0
\end{eqnarray}
has a solution in $K^n$ iff for any $N$ the linear hull of the rows
of $C_N$ does not contain the vector $(1,0,\dots,0)$. An { \it
effective} Nullstellensatz provides an upper bound on $N$ for which
the latter equivalence holds. The bound $N< (\max_j
\{\deg(F_j)\})^{O(n)}$ close to optimal was obtained in
\cite{Heintz}, \cite{Kollar}.

Thus, the effective Nullstellensatz is equivalent to the following.
System (\ref{1}) has a solution iff the linear system $C_N\cdot
(y_1,y_2,\dots)=0$ has a solution with $y_1\neq 0$ for an
appropriate $N$ depending on $n$ and on $\max_j \{\deg(F_j)\}$. We
call the latter statement the {\it effective dual Nullstellensatz}.
The equivalence that (\ref{1}) has a solution iff the linear system
$C_N\cdot (y_1,y_2,\dots)=0$ has a solution with $y_1\neq 0$ for any
$N$, we call the {\it dual Nullstellensatz}. Finally, the statement
(also equivalent to Nullstellensatz) that (\ref{1}) has a solution
iff the infinite linear system $C\cdot (y_1,y_2,\dots)=0$ has a
solution with $y_1\neq 0$, we call the {\it infinite dual
Nullstellensatz}. The latter infinite linear system makes sense
because each row of $C$ contains just a finite number of non-zero
entries.

\section{Conjecture on a tropical dual
Nullstellensatz}\label{section2}

Below we assume that the tropical semi-ring $T=\RR _{\infty}$, but
for the sake of simplifying the exposition we study tropical zeroes
defined over $\RR$ (although, one could also consider zeroes defined
over $\RR_{\infty}$). For each monomial $Q_l=a_l\otimes X_1^{\otimes
i_{1,l}}\otimes \cdots \otimes X_n^{\otimes i_{n,l}}$ of a tropical
polynomial $f=\bigoplus_l Q_l$ we plot the point
$(a_l,i_{1,l},\dots,i_{n,l})\in \RR\times \ZZ^n \subset \RR^{n+1}$.
Then a point $x=(x_1,\dots,x_n)\in \RR^n$ is a tropical zero of $f$
iff the linear function $(a,i_1,\dots,i_n)\to a+i_1\cdot x_1+\cdots
i_n\cdot x_n$ attains its minimum at the plotted points at least
twice.

Therefore, without changing the set of tropical zeroes of $f$ one
can replace the plotted points by their convex hull. Moreover,
w.l.o.g. for any point $(a,b_1,\dots,b_n)\in \RR^{n+1}$ from this
convex hull one can add the ray $\{(b,b_1,\dots,b_n):\, b\ge a\}$.
The resulting convex set $P(f)\subset \RR^{n+1}$ we call the
(extended) {\it Newton polyhedron} of $f$. Thus, w.l.o.g. one can
modify $f$ replacing it by a tropical polynomial whose plotted
points are just the points of the form $(a,i_1,\dots,i_n)\in
(\RR\times \ZZ^n)\cap P(f)$ with the minimal possible $a$. Finally,
so modified tropical polynomial has the same set of tropical zeroes
as $f$, and (in abuse of notations) we keep for it the same
notation. We say that the modified tropical polynomial is in the
{\it convex form}, and from now on we consider tropical polynomials
only in the convex form. Observe that $x$ is a tropical zero of $f$
iff for the maximal $b\in \RR$ such that the hyperplane
$\{(z_1,\dots,z_{n+1}): \, z_1+x_1\cdot z_2+\cdots+x_n\cdot
z_{n+1}=b\}\subset \RR^{n+1}$ has a non-empty intersection with
$P(f)$, the hyperplane has at least two common points with $P(f)$.

Similarly to the classical algebra to a system of tropical
polynomials
\begin{eqnarray}\label{2}
f_1,\dots,f_s
\end{eqnarray}
in $n$ variables we associate the Cayley matrix
$C:=C(f_1,\dots,f_s)$ over $\RR_{\infty}$ consisting of the
coefficients of (\ref{2}). The columns of $C$ correspond to the
tropical monomials of the form $X^{\otimes I},\, I\in \ZZ^n$, and
the rows of $C$ correspond to the tropical polynomials of the form
$X^{\otimes I}\otimes f_j,\, I\in \ZZ^n,1\le j\le s$. Note that
unlike the classical algebra the tropical Cayley matrix is infinite
in all 4 directions.

\vspace{3mm} {\bf Conjecture 1 on a tropical infinite dual
Nullstellensatz}. System (\ref{2}) has a tropical zero iff the
matrix $C$ has a tropical zero.

\vspace{3mm} The latter statement is obvious in the direction that
if (\ref{2}) has a zero then $C$ has a zero (the similar is true for
two conjectures below as well).

Observe that being a particular case of tropical polynomials (of the
tropical degree 1) matrix $C=(c_{i,I})$ (or in other words, a
tropical linear system) has a tropical zero $(\dots,y_I,\dots)$ if
for every row $i$ of $C$ (in the language of classical algebra) the
minimum $\min_I\{c_{i,I}+y_I\}$ is attained at least for two
different coordinates $I$. Note that a tropical zero of $C$ makes
sense because every row of $C$ contains just a finite number of
finite (so, from $\RR$) entries.

Similarly to the classical algebra for an integer $N$ denote by
$C_N$ a (finite) submatrix of $C$ formed by the rows $X^{\otimes
I}\otimes f_j, \, I=(i_1,\dots,i_n)\in \ZZ^n,\, 1\le j\le s$ with
$|i_1|+\cdots+|i_n|\le N$, and by the columns of $C$ which contain
at least one finite entry at one of these rows.

\vspace{3mm} {\bf Conjecture 2 on a tropical dual Nullstellensatz}.
System (\ref{2}) has a tropical zero iff for any $N$ the matrix
$C_N$ has a tropical zero.

\vspace{3mm} {\bf Conjecture 3 on a tropical effective dual
Nullstellensatz}. There is a function $N$ on $n$ and on $\trdeg
(f_j),\, 1\le j\le s$ such that (\ref{2}) has a tropical zero iff
the matrix $C_N$ has a tropical zero.

\vspace{3mm} Clearly, Conjecture 3 implies Conjecture 2, which in
its turn implies Conjecture 1.

\section{Convex-geometric rephrasing of the tropical dual
Nullstellensatz}\label{section3}

In the present Section we give a rephrasing of Conjecture 1 (and
similarly of Conjectures 2, 3) in terms of the convex geometry in
$\RR^{n+1}$. Thus, assume that Cayley matrix $C$ has a tropical zero
$(\dots,y_I,\dots),\, I\in \ZZ^n$.

For any $I\in \ZZ^n$ consider the shift $P(f_j)+(0,I)\subset
\RR^{n+1},\, 1\le j\le s$ of the Newton polyhedron. We say that a
set $U\subset \RR^{n+1}$ {\it lies above} (with respect to the first
coordinate) a set $V\subset \RR^{n+1}$ if for any pair of points
$(u_1,w_1,\dots,w_n)\in U,\, (v_1,w_1,\dots,w_n)\in V$ we have
$u_1\ge v_1$.

\begin{proposition}\label{proposition}
The following statement is equivalent to Conjecture 1.

For $I\in \ZZ^n,\, 1\le j\le s$ take the minimal $a\in \RR$ such
that the polyhedron $P(f_j)+(a,I)$ lies above the set $Y:=\{(-y_J,\,
J): \, J\in \ZZ^n\}$. Assume that for any $I\in \ZZ^n,\, 1\le j\le
s$ the polyhedron and $Y$ have at least two common points. Then
there exists a hyperplane $H\subset \RR^{n+1}$ defined by a linear
equation $z_1+b_2\cdot z_2+\cdots+b_n\cdot z_n=0$ such that for each
$1\le j\le s$ for the minimal $b\in \RR$ with the property that the
polyhedron $P(f_j)$ lies above the hyperplane $H-(b,0)$, the
intersection of $P(f_j)$ and $H-(b,0)$ has at least two points.
\end{proposition}

For an equivalent statement to Conjecture 2 one has for any $N$ to
consider all $I=(i_1,\dots,i_n)$ such that $|i_1|+\cdots+|i_n|\le
N$. Respectively, for Conjecture 3 one has to take $N$ as a suitable
function in $n$ and in $\trdeg(f_j),\, 1\le j\le s$.

\section{Tropical effective dual Nullstellensatz for univariate
polynomials}\label{section4}

Now let $n=1$. In this case for a pair of tropical polynomials
$f_1,\, f_2$ a tropical effective dual Nullstellensatz follows from
\cite{Tabera} with the bound $N\le \trdeg (f_1)+ \trdeg (f_2)$, but
since this approach relies on the (classical) resultant of a pair of
(classical) polynomials being {\it liftings} of $f_1,\, f_2$,
respectively, the approach fails for overdetermined tropical systems
($s\ge 3$).

\begin{theorem}\label{theorem}
A tropical effective dual Nullstellensatz for univariate tropical
polynomials $f_1,\dots,f_s$ holds with $N\le 4\cdot
(\trdeg(f_1)+\cdots+\trdeg(f_s))$.
\end{theorem}
{\bf Proof}. Fix $1\le j\le s$ for the time being. For the convex
polyhedron $P:=P(f_j)\subset \RR^2$ and $i\in \ZZ$ take the minimal
$a_i\in \RR$ such that the shifted polygon $P_i:=P(f_j)+(a_i,i)$
lies above the set $Y=\{(-y_l,l):\, l\in \ZZ\}$ (see
Proposition~\ref{proposition}). By the assumption for any $i\in \ZZ$
there exist at least two points $(u_1,l_1),\, (u_2,l_2) \in
P_i \cap Y,\, l_1<l_2$. Points from the latter intersection we call
{\it extremal points} of $P_i$.

\begin{lemma}\label{lemma1}
The function $i\to a_i$ is convex.
\end{lemma}

{\bf Proof of Lemma~\ref{lemma1}}. Suppose the contrary and let
$2\cdot a_i>a_{i-1}+a_{i+1}$ for a certain $i$. Let $(u_1,l_1),\,
(u_2,l_2) \in
P_i \cap Y$. Denote by $$S=\{(v,w):\, v-w\cdot (a_i-a_{i-1})\le
u_1-l_1\cdot (a_i-a_{i-1}),\, v+w\cdot (a_i-a_{i+1})\ge u_1+l_1\cdot
(a_i-a_{i+1})\}\subset \RR^2$$ \noindent the sector with the vertex
at the point $(u_1,l_1)$ between two rays
$R_+=(u_1,l_1)+\{\lambda\cdot (a_i-a_{i+1},-1):\, \lambda \ge 0\}$
and $R_-=(u_1,l_1)+\{\lambda\cdot (a_i-a_{i-1},1):\, \lambda \ge
0\}$. We claim that
$P_i\subset S$.

Indeed, consider a left adjacent to $(u_1,l_1)$ point
$(u_+,l_1-1)\in \partial P_i$
 on the boundary of $P_i$
(provided that such a point does exist). If $u_+<u_1+l_1\cdot
(a_i-a_{i+1})$ (in other words, the point $(u_+,l_1-1)$ lies
strictly below the ray $R_+$, cf. the description of $S$) then the
point $(u_+,l_1-1)+(a_{i+1}-a_i,1)\in P_{i+1}$
lies strictly below $Y$, the achieved contradiction implies that
$(u_+,l_1-1)\in S$. In a similar way a right adjacent to $(u_1,l_1)$
point $(u_-,l_1+1)\in \partial P_i$
(provided that it
does exist) belongs to $S$, which justifies the claim.

By the same token the parallel shift $S+(u_2,l_2)-(u_1,l_1)$ of the
sector $S$ (with its vertex at the point $(u_2,l_2)$) also contains
$P_i$. This contradicts to the convexity of
$P_i$ and completes the proof of Lemma~\ref{lemma1}. \bull

\vspace{3mm} Denote by $E:=E(f_j)\subset \RR^2$ the polygon with
the vertices in the extremal points of
$P_i$ for all $i\in \ZZ$. Below we enumerate the (finite) edges of
the polygon
$P$ from the left to the right. Denote by $(b_r,1)$ the vector
parallel to the $r$-th edge of $P$.

\begin{lemma}\label{lemma2}
Let $(u_1,l_1),\dots,(u_t,l_t)\in
 P_i,\, l_1<\cdots <
l_t$ be all the extremal points of
$P_i$. Let
the point $(u_t,l_t)-(a_i,i)\in
P$ lie in the $r$-th (finite) edge of
$P$ (when the latter point belongs to the $r$-th and to the
$(r+1)$-th edges we agree that the point lies in the $r$-th edge).
Then $a_{i+1}-a_i\ge b_r$.

For any extremal point $(v,k)$ of
$P_{i+1}$ the point $(v,k)-(a_{i+1},i+1)\in
P$ lying in the $q$-th edge of
$P$ either satisfies an inequality $q\ge r$ or $(v,k)-(a_{i+1},i+1)$
is the vertex of the $(r-1)$-th and $r$-th edges of
$P$ (in the latter case $(v,k)$ is the leftmost extremal point of
$P_{i+1}$). There exists an extremal point $(v,k)$ for which either
$q=r$ and $(v,k)-(a_{i+1},i+1)$ not being the vertex of the $r$-th and
$(r+1)$-th edges of the polygon
$P$ or $(v,k)-(a_{i+1},i+1)$ is the vertex of the $(r-1)$-th and $r$-th
edges of
$P$ iff $a_{i+1}-a_i= b_r$. Moreover, when $a_{i+1}-a_i= b_r$ any
extremal point $(u_m,l_m)$ of
$P_i$ with $(u_m,l_m)-(a_i,i)$ lying in the $r$-th edge of
$P$ is also an extremal point of
$P_{i+1}$.
\end{lemma}

{\bf Proof of Lemma~\ref{lemma2}}. Consider the point
$(u_t,l_t)-(b_r,1)\in
P_i$. Then the point $((u_t,l_t)-(b_r,1))+(a_{i+1}-a_i,1)\in
P_{i+1}$ should lie above the extremal point $(u_t,l_t)$, this
entails the inequality $a_{i+1}-a_i\ge b_r$.

Let $(v,k)$ be an extremal point of
$P_{i+1}$ with $(v,k)-(a_{i+1},i+1)$ lying in the $q$-th edge of
$P$. The point $(v,k)-(a_{i+1}-a_i,1)$ lies in the $q$-th edge of
the polygon
$P_i$. If $q<r$ and the point $(v,k)-(a_{i+1},i+1)$ is not the
vertex of the $(r-1)$-th and $r$-th edges of
$P$ then its shift $(v,k)=((v,k)-(a_{i+1}-a_i,1))+(a_{i+1}-a_i,1)$
lies {\it strictly} inside the polygon
$P_i$, and therefore $(v,k)$ can not be an extremal point. The
achieved contradiction implies that either $q\ge r$ or
$(v,k)-(a_{i+1},i+1)$ is the vertex of $(r-1)$-th and $r$-th edges
of
$P$.

 When $a_{i+1}-a_i > b_r$ a similar argument shows that either
$q>r$ or $(v,k)-(a_{i+1},i+1)$ is the vertex of the $r$-th and
$(r+1)$-th edges of
 $P$. Finally, when $a_{i+1}-a_i=b_r$, for
any extremal point $(u_m,l_m)$ of
$P_i$ with $(u_m,l_m)-(a_i,i)$ lying in the $r$-th edge of
$P$ take the point $(u_m,l_m)-(a_{i+1}-a_i,1)\in
P_i$, then the point
$(u_m,l_m)=((u_m,l_m)-(a_{i+1}-a_i,1))+(a_{i+1}-a_i,1)\in
 P_{i+1}$ is also an extremal point of
$P_{i+1}$. \bull

\begin{remark}
Lemma~\ref{lemma2} is formulated for the shifts passing from the
polygon
$P_i$ to
$P_{i+1}$ (so, from the left to the right). By the same token a
similar statement holds while passing from
$P_{i+1}$ to
$P_i$  (so, from the right to the left).
\end{remark}

\begin{lemma}\label{lemma3}
The polygon $E$ is convex.
\end{lemma}

{\bf Proof of Lemma~\ref{lemma3}}. We prove by induction on $i$ that
the union of the extremal points of $P_0,\dots,P_i$ form a convex
polygon $E_i$. At the inductive step we consider the polygon
$P_{i+1}$ (so, we move from the left to the right). By the same
token one can alternatively consider the polygon $P_{-1}$ (so, move
from the right to the left). This would entail Lemma~\ref{lemma3}.

We say that a polygon with the vertices
$\dots,w_{i-1},w_i,w_{i+1},\dots$ is convex at the vertex $w_i$ if
there is a line passing through $w_i$ such that both points
$w_{i-1},\, w_{i+1}$ lie above this line.

Let $(u_1',l_1'),\, (u_2',l_2'),\dots,\, l_1'<l_2'<\cdots$ be all
the extremal points of $P_{i+1}$ (if they exist) being not extremal
points of $P_i$. Lemma~\ref{lemma2} implies that the point
$(u_1',l_1')-(a_{i+1},i+1)$ either lies in the $q$-th edge of $P$
with $q>r$ or it is the vertex of the $r$-th and $(r+1)$-th edges of
$P$ (we keep the notations from Lemma~\ref{lemma2}).

The inductive hypothesis and Lemma~\ref{lemma2} entail that the
polygon $E_{i+1}$ is convex at all the vertices of $E_i$, perhaps,
with the exception of the rightmost extremal point $(u_t,l_t)$ of
$E_i$ (and simultaneously of $P_i$). The point $(u_t,l_t)$ lies in
the $r$-th edge of the polygon $P_i$, and both polygons $P_i,\,
P_{i+1}$ lie above the line $L$ spanned by this edge (due to
Lemma~\ref{lemma2}), whence $E_{i+1}$ is convex at its vertex
$(u_t,l_t)$ as well.

Since the extremal points $(u_1',l_1'),\, (u_2',l_2'),\dots$ are
located on the convex polygon $P_{i+1}$ we get that $E_{i+1}$ is
convex at its vertices $(u_2',l_2'),\dots$. Thus, it remains to
verify that $E_{i+1}$ is convex at its vertex $(u_1',l_1')$.

Denote the vector $w:=(u_2',l_2')-(u_1',l_1')$. The points
$p:=(u_1',l_1')-(a_{i+1}-a_i,1),\, (u_2',l_2')-(a_{i+1}-a_i,1)\in
P_i$. Therefore, the point $p$ lies in a sector $S_0$ with the
vertex $(u_t,l_t)$ formed by the rays $(u_t,l_t)+\{\lambda \cdot
(b_r,1):\, \lambda \ge 0\}\subset L$ and $(u_t,l_t)+\{\lambda \cdot
w:\, \lambda \ge 0\}$. Now consider a sector $S_1\subset S_0$
parallel to $S_0$ with the vertex $p$ formed by the rays
$p+\{\lambda \cdot (b_r,1):\, \lambda \ge 0\}$ and $p+\{\lambda
\cdot w:\, \lambda \ge 0\}$. The point
$(u_1',l_1')=p+(a_{i+1}-a_i,1)$ is located in $S_1$ due to
Lemma~\ref{lemma2} and taking into the account that the point
$(u_1',l_1')$ is extremal in $P_{i+1}$ and thereby, can not lie
strictly inside $P_i$. Hence the polygon $E_{i+1}$ is convex at its
vertex $(u_1',l_1')$. \bull

\begin{remark}\label{remark2}
The latter statement that $E_{i+1}$ is convex at its vertex
$(u_1',l_1')$ becomes obvious when the point $(u_t,l_t)$ is an
extremal point of $P_{i+1}$, this is equivalent to the equality
$a_{i+1}-a_i=b_r$ due to Lemma~\ref{lemma2}. In case when
$a_{i+1}-a_i>b_r$ the polygon $P_{i+1}$ has no common extremal
points with $E_i$.
\end{remark}

\begin{corollary}\label{corollary1}
Any edge $e=((u,l),\, (u',l'))$ of the convex polygon $E$ is one of
the following three types:

1) either $(u,l),\, (u',l')\in P_i$ for a certain $i\in \ZZ$ where
the point $(u,l)$ lies in the $r$-th edge of $P_i$, the point
$(u',l')$ lies in the $r'$-th edge of $P_i$ for some $r<r'$, except
the case when $(u,l)$ is the vertex of the $(r-1)$-th and $r$-th edges
of $P_i$ and $(u',l')$ lies in the $r$-th edge of $P_i$ (in the latter
case $e$ is parallel to the $r$-th edge of $P_i$, cf. 3) below);

2) either the point $(u,l)$ lies in the $r$-th edge of $P_i$ for a
certain $i\in \ZZ$, the point $(u',l')$ lies in the $r'$-th edge of
$P_{i+1}$ for some $r,r'$, and the point $(u',l')-(a_{i+1}-a_i,1)$
lies in the $r'$-th edge of $P_i$, moreover either $r<r'$ or
$(u',l')-(a_{i+1}-a_i,1)$ is the vertex of the $r$-th and
$(r+1)$-th edges of $P_i$. Case 2) occurs when $a_{i+1}-a_i>b_r$
(see Lemma~\ref{lemma2} and Remark~\ref{remark2});

3) or $e$ is parallel to an edge of $P$.
\end{corollary}

The edges $e$ of $E$ of types 1), 2) we call {\it intermediate} and
the edges of types 3) we call $r$-{\it principal} when $e$ is
parallel to the $r$-th edge of $P$. For an edge of the type either
1) or 2) we define its {\it projection} (to the second coordinate)
as the interval $(l-i,l'-i)$ for the type 1) and as $(l-i,l'-i-1)$
for the type 2).

\begin{corollary}\label{corollary2}
i) For a pair of adjacent intermediate edges of $E$ their
projections are also adjacent (in the same order).

ii) For each $r$ all $r$-principal edges of $E$ (if they exist)
constitute an interval in $E$ (parallel to the $r$-th edge of $P$).
We call it $r$-{\it interval}. Among these intervals there are
either two intervals infinite in one of directions or one interval
infinite in both directions.

Let $e_-:=((u_-,l_-),\, (u_-',l_-'))$ be an edge of $E$ adjacent to
the $r$-interval from the left (provided that the $r$-interval is
not infinite to the left). Then $e_-$ is intermediate. Assume
for definiteness that $(u_-',l_-')\in P_i$ for a certain $i$, while
either $(u_-,l_-)\in P_i$ in case of the type 1) (see
Corollary~\ref{corollary1}) or $(u_-,l_-)\in P_{i-1}$ in case of the
type 2). Then the point $(u_-',l_-')$ lies in the $r$-th edge of
$P_i$.

Similarly, let an edge $e_+:=((u_+,l_+),\, (u_+',l_+'))$ be an edge
adjacent to the $r$-interval from the right (provided that the
$r$-interval is not infinite to the right). Then $e_+$ is
intermediate. Assume that $(u_+,l_+)\in P_i$ for a certain $i$ and
either $(u_+',l_+')\in P_i$ in case of the type 1) or
$(u_+',l_+')\in P_{i+1}$ in case of the type 2). Then the point
$(u_+,l_+)$ either lies in the $r$-th edge of $P_i$ or $(u_+,l_+)$
is the vertex of the $(r-1)$-th and $r$-th edges of $P_i$.

iii) Denote by $(d_1,k_1),\, (d_2,k_2)$ the endpoints of the $r$-th
edge of $P$. Then for any pair of adjacent extremal points
$(d_1',k_1'),\, (d_2',k_2')$ in the $r$-interval of $E$ we have
$k_2'-k_1'\le k_2-k_1$.
\end{corollary}

Finally, we complete the proof of Theorem~\ref{theorem}. So far, we
studied the convex polygon $E(f_j)$ for a fixed $1\le j\le s$ (see Lemma~\ref{lemma3}).
Now we consider the intersection ${\cal E}:= \bigcap _{1\le j\le s}
E(f_j)$. Every edge $\epsilon$ of the convex polygon $\cal E$ is some subinterval of
either an intermediate edge of $E(f_j)$ or an $r$-interval for
certain $1\le j\le s$ and $r$. The total sum of the lengths of the
projections of the edges being subintervals of intermediate edges of
$E(f_j),\, 1\le j\le s$ does not exceed $3\cdot \sum_{1\le j\le s}
\trdeg f_j$ due to i), ii) of Corollary~\ref{corollary2}.

Observe that if $\epsilon$ is a subinterval of an $r$-interval of
$E(f_j)$ for a certain $1\le j\le s$ and not all the points of
$\epsilon$ belong to all polygons $E(f_{j_1}),\, 1\le j_1\le s$ (the
latter is equivalent to that any strictly inside point of $\epsilon$
does not belong to all $E(f_{j_1}),\, 1\le j_1\le s$) then
$\epsilon$ can not contain extremal points strictly inside itself.
Hence the total sum of the lengths of the projections of all the
edges of $\cal E$ being subintervals of some $r$-intervals of
$E(f_j),\, 1\le j\le s$ does not exceed $\sum_{1\le j\le s} \trdeg
f_j$ by virtue of iii) of Corollary~\ref{corollary2}.

Thus, a truncation ${\cal E}_N$ of $\cal E$ with the length of the
projection to the second coordinate equal $N$, where $N\ge 4\cdot
\sum_{1\le j\le s} \trdeg f_j$, contains an edge which is a common
subinterval of $r_j$-intervals for appropriate $r_j$ of all
$E(f_j),\, 1\le j\le s$. Taking into the account the
Proposition~\ref{proposition} we conclude with Theorem~\ref{theorem}.
\bull

\vspace{3mm} It would be interesting to improve the factor $4$ in
Theorem~\ref{theorem}.

 { \vspace{3mm} {\bf Acknowledgements}. The author is grateful
to the Max-Planck Institut f\"ur Mathematik, Bonn for its
hospitality during writing this paper.

\end{document}